%% file: Penult.tex
\documentclass[12pt,a4paper,oneside]{article}
\usepackage[utf8]{inputenc}
\usepackage[english]{babel}
\usepackage[thmmarks]{ntheorem}
\usepackage{amsmath}
\usepackage{amssymb}
\usepackage{graphicx,xcolor}
\usepackage[font=small,labelsep=endash]{caption}
\usepackage[round]{natbib}
\usepackage{calc}
\usepackage[linkcolor=black, citecolor=black,%
unicode, bookmarks=false, hypertexnames=false]{hyperref}

\theoremstyle{break}
\theoremsymbol{}
\newtheorem{theorem}{Theorem}[section]

\theorembodyfont{\upshape}
\theoremstyle{plain}
\renewcommand*{\qedsymbol}{$\blacksquare$}
\theoremsymbol{\qedsymbol}
\newtheorem*{proof*}{Proof}

\input macros.tex

\renewcommand{\textwidth}{6.3in}
\renewcommand{\textheight}{21cm}
\newcommand{\gras}{\boldsymbol}
\numberwithin{equation}{section}
\begin{document}
\title{Penultimate analysis of the conditional multivariate extremes tail model}
\author{T. Lugrin, A. C. Davison, J. A. Tawn}
\date{\today}
\maketitle
\begin{abstract}
  Models for extreme values are generally derived from limit results, which are meant to be good enough approximations when applied to finite samples.
  Depending on the speed of convergence of the process underlying the data, these approximations may fail to represent subasymptotic features present in the data, and thus may introduce bias.
  The case of univariate maxima has been widely explored in the literature, a prominent example being the slow convergence to their Gumbel limit of Gaussian maxima, which are better approximated by a negative Weibull distribution at finite levels.
  In the context of subasymptotic multivariate extremes, research has only dealt with specific cases related to componentwise maxima and multivariate regular variation.
  This paper explores the conditional extremes model \citep{HeffernanTawn2004} in order to shed light on its finite-sample behaviour and to reduce the bias of extrapolations beyond the range of the available data.
  We identify second-order features for different types of conditional copulas, and obtain results that echo those from the univariate context. These results suggest possible extensions of the conditional tail model, which will enable it to be fitted at less extreme thresholds.
  
\end{abstract}

\section{Introduction}
Large-scale catastrophic events can have a major impact on physical infrastructure and society. Multivariate extreme value models help to capture the structure of such events and are used to extrapolate measures of combined risk beyond the range of the available data. Applications include river flooding \citep{KatzParlangeNaveau2002,KeefSvenssonTawn2009,KeefTawnSvensson2009,AsadiDavisonEngelke2015}, extreme rainfall \citep{ColesTawn1996JRSSB,SuvegesDavison2012,HuserDavison2014}, wave height and extreme sea surge \citep{deHaandeRonde1998} and high concentrations of air pollutants \citep{HeffernanTawn2004,EastoeTawn2009}. Many applications have also contributed to the improvement of risk assessment in finance \citep{PoonRockingerTawn2003,HilalPoonTawn2011,HilalPoonTawn2014}. 

All these approaches can be characterised in terms of their ability to model probabilities in the joint tail. Specifically, 
for random variables $X_1$ and $X_2$ with marginal distributions $F_1$ and $F_2$, a basic summary of extremal dependence is \citep{ColesHeffernanTawn1999} 
\begin{equation}\label{eq:chi}
  \chi = \lim_{p\to 1}\pr\left\{F_2(X_2) >  p \mid F_1(X_1) > p\right\}.
\end{equation}  
In the asymptotic dependence case,  $\chi > 0$ and the largest values of $X_1$ and $X_2$ can occur together, whereas in the asymptotic independence case, $\chi = 0$ and these largest values cannot  occur simultaneously.
Many extreme value models 
can only capture asymptotic dependence, but this excludes important cases of asymptotic independence; for example, $\chi=0$ for all Gaussian copulas with correlation $\rho<1$.
In the case when $\chi = 0$, the key to determining the form of the extremal dependence is the rate at which the conditional probability in~\eqref{eq:chi} tends to zero, given by \citet{ColesHeffernanTawn1999} as a complementary extremal dependence measure $\bar\chi\in(-1 ,1]$, which for Gaussian copulas gives $\bar\chi = \rho$.

Modelling approaches proposed by \citet{LedfordTawn1997}, \citet{HeffernanTawn2004}, \citet{WadsworthTawnDavisonElton2017} and \citet{HuserWadsworth2018} cover both asymptotic dependence and asymptotic independence. In all cases, an asymptotic form for the joint tail is specified above a high threshold and for inference and extrapolation of extreme probabilities the approach relies on the assumption that the model holds exactly. To illustrate the implications of this assumption, consider the case of asymptotic dependence, in which $\pr\{F_2(X_2) >  p \mid F_1(X_1) > p\}\equiv\chi$ for all $p$ above the selected threshold once the model and threshold are selected, so there is no scope to allow for a penultimate value for these probabilities that converges to $\chi$ as $p\rightarrow 1$. Similar issues arise when estimating $\bar\chi$.

The focus of this paper is to try to establish penultimate models that are appropriate in the most flexible existing model for multivariate extremes, that of conditional extremes introduced by \citet{HeffernanTawn2004}.

In our analysis of conditional extremes, we focus on the bivariate case to simplify the notation; extension to the general multivariate case is straightforward  \citep{HeffernanTawn2004}. The conditional model was originally presented for marginally Gumbel distributed random vectors, but \citet{KeefPapaTawn2013} showed that formulation on the Laplace scale is simpler when positive or negative dependence is possible. 
Thus we first transform our variables $(X_1,X_2)$ to random variables $(X,Y)$ with Laplace margins via the probability integral transform,
\[
  X = \sign\{1-2F_1(X_1)\}\log \left[1 - \left\{\left\vert 1 - 2F_1(X_1)\right\vert\right\}\right],
\]
and similarly for $Y$, preserving the dependence structure through the copula, according to Sklar's~(\citeyear{Sklar1959}) representation theorem.

The \citet{HeffernanTawn2004} approach presupposes that there exist normalising functions $a(\cdot)$ and $b(\cdot)>0$ such that
\begin{equation}\label{eq:HT_random_norming}
    \pr\left\{\dfrac{Y-a(X)}{b(X)}
         \leq z,\, X - u > x \mmid X > u\right\} \to H(z)\exp(-x),\quad u\to\infty,
\end{equation}
where $H(\cdot)$ is a non-degenerate distribution function with no mass at infinity.
Under mild assumptions on the joint distribution of $(X,Y)$, \citet{HeffernanResnick2007} show that the normalising functions are of the form
\begin{equation}\label{eq:slowly-varying_normings}
  a(x) = x {\cal L}_a(x),\quad b(x) = x^\beta {\cal L}_b(x),\ \beta < 1,
\end{equation}
with ${\cal L}_a(x)$ and ${\cal L}_b(x)$ two slowly-varying functions, i.e., satisfying ${\cal L}(xt)/{\cal L}(x)\to 1$ for any  $t > 0$ as $x\to\infty$.
By considering a broad class of parametric copula  models, \citeauthor{HeffernanTawn2004} derive parametric forms for $a(\cdot)$ and $b(\cdot)$ that yield a parsimonious model and cover a broad range of extremal dependence structures not described by models arising from the standard theory for multivariate extremes, namely
\begin{equation}\label{eq:ultimate_normings}
\begin{aligned}
  a_0(x) &= \alpha x,\quad \alpha\in [-1,1],\\
  b_0(x) &= x^\beta,\quad \beta\in (-\infty,1).
\end{aligned}
\end{equation}
This corresponds to approximating the slowly-varying functions by ${\cal L}_a(x) \equiv \alpha$ and ${\cal L}_b(x) \equiv 1$. Setting ${\cal L}_b(x) = 1$ is equivalent to setting ${\cal L}_b(x)=b$ for any  positive constant $b$, with the change in norming absorbed into the variance of $H(\cdot)$. If $\alpha = 1$ and $\beta = 0$, then $(X,Y)$  are asymptotically dependent with $\chi = \int_0^\infty \bar H(-z)e^{-z}dz$, and otherwise they are asymptotically independent.
To form a statistical model, \citet{HeffernanTawn2004}  assume that limit~\eqref{eq:HT_random_norming} holds exactly above some finite $u$ with  norming functions of the form~\eqref{eq:ultimate_normings}.

\citet{PapaTawn2016} found inverted max-stable processes for which the norming functions of the form~\eqref{eq:ultimate_normings} are inadequate and more general functions ${\cal L}_a(x)$ and ${\cal L}_b(x)$ of the form \eqref{eq:slowly-varying_normings} are required. It is natural therefore to question whether we can extend the functions $a_0(x)$ and $b_0(x)$
to give better finite approximations.

Our main contribution is to derive the subasymptotic behaviour of the conditional tail approach for three copulas that span diverse extremal dependence structures.
Our objective is to identify a second-order behaviour across conditional copulae from which we can derive a general penultimate conditional model for extremes.
The core value of the work is to suggest in Section~\ref{sec:discussion} new penultimate
forms for $a_0(\cdot)$ and $b_0(\cdot)$ that can be used to broaden the simple parametric family~\eqref{eq:ultimate_normings}, and help reduce extrapolation bias. This is particularly important as small differences in the parameter estimates at finite levels can result in large differences in extreme risk measures.

Another aspect of the subasymptotic behaviour is the limiting conditional independence of $Z \sim H$ with $X - u$, given $X > u$, as $u \to \infty$. At a subasymptotic level, the distribution of $X - u$ is exponential, but that of $Z$ depends on $X$. This subasymptotic distribution, $H_x(\cdot)$, say, tends to $H(\cdot)$ as $x \to \infty$, so the independence property is lost in this penultimate model which might therefore enable it to provide a better fit at a lower threshold than limit models.

Such subasymptotic behaviour is also required when conducting simulation studies in order to assess the performance of methods to fit the conditional model, as the estimates of $a_0(x)$ and $b_0(x)$ for $x>u$  can misleadingly suggest a poor fit
as it is $a(x)$ and $b(x)$ for $x>u$ that are being estimated.

Although the focus of this paper is the extremal dependence structure, we outline the established convergence of the margins to their limiting distributions in Section~\ref{sec:existing} to set a framework for our extremal dependence analysis. The rest of the paper is organised as follows.
In Section~\ref{sec:existing}, we review penultimate analyses in bivariate contexts.
Section~\ref{sec:conditional_extremes} introduces the framework used to study the penultimate behaviour of the conditional tail model.
In Section~\ref{sec:examples}, we consider three copulas, with proofs of results given in the Appendix.
In Section~\ref{sec:discussion}, we summarise our findings through penultimate parametric models which extend the Heffernan--Tawn class of norming functions.

\section{Existing penultimate analyses}\label{sec:existing}

The founding theorem in the theory of univariate extreme values  characterises the limiting distribution $G_\xi(x)$ of  the maximum of a series of independent and identically distributed univariate random variables $X_1,\ldots,X_n\sim F$,  suitably normalised by sequences $a_n > 0$ and $b_n$, namely, as $n\to\infty$,
\[
  \pr\left(\dfrac{\max\left\{X_1,\ldots,X_n\right\} - b_n}{a_n}\leq x\right)
  = F^n\left(a_nx + b_n\right)
  \to G_\xi(x)
  = \exp\left\{-(1+\xi x)_+^{-1/\xi}\right\},
\]
with $a_+ = \max(a,0)$.
In practice the limit distribution $G_\xi(x)$ is assumed to be exact for large $n$, by taking $F^n(x) = G_\xi\{(x-b_n)/a_n\}$, with $a_n$, $b_n$ and $\xi$ parameters to be estimated. This leads to the generalised extreme value distribution being proposed as a parametric model for inference and extrapolation based on block maxima \citep{ISMEV}.

\citet{FisherTippett1928} raised the question of the accuracy of this approximation in the Gaussian case, i.e., $F(x)=\Phi(x)$.  Then $G_\xi(x)$ is Gumbel, i.e., $\xi = 0$, but  they showed  that at finite levels $\Phi^n(a_nx+b_n)$ is better approximated by a generalised extreme value distribution with shape parameter $\xi<0$, a negative Weibull distribution.

The approximation of $F^n(a_nx+b_n)$ by $G_\xi(x)$ in extreme value applications is of concern when the convergence of the former to the latter is particularly slow, as any inaccuracy in the parameter estimates amplifies and introduces bias in extrapolations.

The study of rates of convergence of $\vert F^n(a_nx+b_n)-G_\xi(x)\vert$ towards zero has a long history. The first attempts to characterise this for any value of $\xi$ date to \citet{Gomes1984,Gomes1994} and unpublished work of~\citet{Smith1987}. 
Assuming the existence of the density $f(x)=F'(x)$ and its derivative $f'(x)$, we define the reciprocal hazard function $h(x)=\{1-F(x)\}/f(x)$, and we know from the \citet{vonMises1936} conditions that
\begin{equation}\label{eq:ultimate_xi}
  \xi = \lim_{x\to x^F} h'(x),
\end{equation}
where $x^F$ is the upper support point of $f$.
From \eqref{eq:ultimate_xi}, it is natural to consider the sub-asymptotic shape parameter $\xi_n = h'(b_n)$, with $F(b_n)=1-1/n$. With this sub-asymptotic $\xi_n$, \citet{GomesPestana1987} and \citet{Gomes1994} show that
\begin{equation}\label{eq:univariate_convergence}
  F^n(a_nx+b_n)-G_\xi(x) = O(\xi_n-\xi),\quad n\to\infty,
\end{equation}
and give the structure of the remainder term on the right-hand side for a broad class of distribution functions $F$, including $\Phi$ \citep{Anderson1971}. For this  class, \citet{Gomes1994} also gives
\begin{equation}\label{eq:univariate_subasymptotic_convergence}
  F^n(a_nx+b_n) - G_{\xi_n}(x) = O\left\{(\xi_n-\xi)^2\right\},\quad n\to\infty,
\end{equation}
so replacing $\xi$ by $\xi_n$ can greatly improve the rate of convergence.

To illustrate~\eqref{eq:univariate_convergence} and \eqref{eq:univariate_subasymptotic_convergence}, consider the standard Gaussian case $F(x)=\Phi(x)$ and $f(x)=\varphi(x)=\Phi'(x)$. Mills' ratio can be used to derive the approximation to the survival distribution function $\bar F(x) = \varphi(x)\left\{x^{-1}+O\left(x^{-3}\right)\right\}$, from which we get the approximate reciprocal hazard function $h(x)=\{1-\Phi(x)\}/\varphi(x) = x^{-1}+O\left(x^{-3}\right)$ and its derivative $h'(x)= -x^{-2}+O\left(x^{-4}\right)$. We also have $b_n=\sqrt{2\log n}+o(1)$  \citep[p.~14]{LeadbetterLindgrenRootzen1983}, from which we can conclude that
\begin{equation}\label{eq:xi_gaussian}
  \xi_n=h'(b_n)=-1/(2\log n)+O\{1/(\log n)^2\}.
\end{equation}
The convergence rate \eqref{eq:univariate_convergence} is $O(1/\log n)$ using the limit shape parameter,  improving to $O\{1/(\log n)^2\}$ when $\xi$ is replaced by $\xi_n$, as in \eqref{eq:univariate_subasymptotic_convergence}. Here the penultimate shape parameter $\xi_n<0$ for finite $n$, in agreement with the observation, due to \citet{FisherTippett1928},  that a negative Weibull-type distribution yields better approximations at finite levels.

The first study that has discussed the penultimate properties of bivariate maxima was \citet{BofingerBofinger1965}, who derive the penultimate correlation of componentwise maxima in bivariate Gaussian samples of sizes $n=2,\ldots,50$; the limit correlation was shown to be zero by \citet{Sibuya1960}, corresponding to $\chi = 0$ in \eqref{eq:chi}.  \citet{Bofinger1970} extends this analysis to the bivariate gamma and  \citet{Morgenstern1956} distributions, and  sheds light on the form of the penultimate correlation for small $n$. 
This work has been further explained for identically distributed random variables $X$ and $Y$ by \citet{LedfordTawn1996}, who derive a model for joint tails of these variables using penultimate properties of $\pr(X>z,Y>z)$ as $z$ tends to the marginal upper endpoint; this results in a  model that smoothly connects perfect dependence and complete independence. In exponential margins, \citet{WadsworthTawn2013} and \citet{deValk2016} extended this model to allow for different decay rates along different rays emanating from the origin, with these rays corresponding to power relationships in Pareto margins.
There is also very strong parallels with work on inference for the spectral measure of multivariate regular variation with second-order features influencing estimates \citep{Resnick2007,CaiEinmahldeHaan2011}.
Despite these developments, there are currently no penultimate results for conditional multivariate extremes to parallel the univariate penultimate theory.

\section{Sub-asymptotic conditional extremes}\label{sec:conditional_extremes}
The conditional limit \eqref{eq:HT_random_norming} encapsulates the limit conditional independence of $Z=\{(Y-a(X))/b(X)$ and the excesses $X-u$  for large  $X$. We first focus on the marginal limiting behaviour of $Z$. According to \citet{HeffernanResnick2007}, \citet{ResnickZeber2014} and \citet{WadsworthTawnDavisonElton2017}, the conditioning in \eqref{eq:HT_random_norming} can be modified to 
\begin{equation}\label{eq:HT_assumption}
  \lim_{x\to\infty}\pr\left\{\dfrac{Y-a(x)}{b(x)}\leq z\mmid X=x\right\}=H(z),
\end{equation}
where $a(x)$, $b(x)$ and the distribution $H(\cdot)$ are the same as in \eqref{eq:HT_random_norming}.

In our penultimate analysis, our goal is to characterise the behaviour of the remainder terms, defined by
\[
  a(x)-a_0(x)\sim r_a(x),\quad   b(x)-b_0(x)\sim r_b(x),\qquad x\to\infty,
\]
using the notation of \eqref{eq:ultimate_normings}.
Specifically, we consider the second-order normalisation for $a(\cdot)$ and $b(\cdot)$, with
\begin{equation}\label{eq:penultimate_normings}
\begin{aligned}
  a_1(x)&=a_0(x)+r_a(x),\\
  b_1(x)&=b_0(x)+r_b(x).
\end{aligned}
\end{equation}
With these penultimate forms, we are able to refine the normalisation of $Y$ in \eqref{eq:HT_assumption}, yielding the subasymptotic conditional distribution
\begin{equation}\label{eq:penultimate_residual_distribution}
  \pr\left\{\dfrac{Y-a_1(X)}{b_1(X)}\leq z\mmid X=x\right\} = H_x(z),\quad x>u,
\end{equation}
with $H_x(z)\to H(z)$ as $u\to\infty$.
 
\citet{HeffernanTawn2004} give the rate of convergence of the conditional distribution for data arising from various copula models in terms of the order of convergence towards zero, as $x\to\infty$, of
\begin{equation}\label{eq:HT_convergence}
  \pr\left\{\dfrac{Y-a_0(X)}{b_0(X)}\leq z\mmid X=x\right\} - H(z),
\end{equation}
with $(X,Y)$ on the Gumbel scale. We consider how much we can improve this when using the penultimate norming, by studying the rate of convergence to zero of
\begin{equation}\label{eq:HT_convergence_penultimate_norming}
  \pr\left\{\dfrac{Y-a_1(X)}{b_1(X)}\leq z\mmid X=x\right\} - H(z).
\end{equation}
We also want to quantify the subasymptotic remainder, using
\begin{equation}\label{eq:HT_remainder_penultimate}
  \sup_{x>u}\left\vert\pr\left\{\dfrac{Y-a_1(x)}{b_1(x)}\leq z\mmid X=x\right\}-H_x(z)\right\vert,
\end{equation}
along the lines of \eqref{eq:univariate_subasymptotic_convergence} in the univariate context.
In all cases we will present these rates on a scale that is invariant to the marginal choice, by converting to a return period $n$, where $\pr(X>x)=n^{-1}$.

\section{Examples}\label{sec:examples}

We first consider the Gaussian copula, with correlation parameter $\rho$, which has $\alpha = \sign(\rho)\rho^2$ and $\beta = 1/2$, for $\rho \neq 0$, for which the convergence towards the limit \eqref{eq:HT_random_norming} was reported by \citeauthor{HeffernanTawn2004} to be the slowest in the examples they considered, namely $O\{\log(\log n)/(\log n)^{1/2}\}$.
Second, we consider another example with asymptotic independence; the inverted logistic dependence structure (a special case in the class of inverted max-stable distributions studied by \citet{PapaTawn2016}). This distribution has $\alpha = 0$ and $\beta = 1 - \gamma$ where $0<\gamma\le 1$ is the logistic parameter with $\gamma = 1$ corresponding to independence.
This distribution has a faster convergence rate than the Gaussian copula, and we shall see that the subasymptotic $H_x$ can have finite support depending on the precise value of the dependence parameter of this copula. \citeauthor{HeffernanTawn2004} reported the rate to be $O(1/\log n)$ in this case.
Third, the max-stable copula with logistic dependence structure represents the case of asymptotic dependence, $\alpha = 1$ and $\beta = 0$, and is a situation where convergence is $O(1/n)$.
Our penultimate forms for $\alpha_1(x)$ and $\beta_1(x)$ reflect these different rates of convergence, with major, minor and no changes found relative to $a_0(x)$ and $b_0(x)$ respectively.

\subsection{Gaussian distribution}
Let $(V,W)$ have a bivariate standard normal distribution with correlation $\rho\neq 0$ and let $(X,Y)$ be its marginal transform to the Laplace scale, 
\[
X =
\begin{cases}
    -\log 2\{1-\Phi(V)\},& V>0,\\
    \log 2\Phi(V),& V\leq 0,
\end{cases}
\]
and similarly for $Y$ as a function of $W$. The dependence structure of $(X,Y)$ is an example where $\chi=0$, provided $\rho<1$, and $\bar\chi=\rho$, $\rho\in[-1,1]$.
  
\begin{theorem}\label{thm:gaussian}
For $(X,Y)$ with a Gaussian dependence structure, we have the ultimate and penultimate normings \eqref{eq:ultimate_normings} and \eqref{eq:penultimate_normings} for $Y\mid X=x$, with $x$ large, are
  \begin{equation}\label{eq:gaussian_penultimate_normings}
    \begin{aligned}
      a_0(x) &= \sign(\rho)\rho^2 x,&
      b_0(x) &= x^{1/2},\\
      a_1(x) &= \sign(\rho)\rho^2 x + \dfrac{(1-\rho^2)}{2}\log x,&
      b_1(x) &= x^{1/2-1/(4x)},
    \end{aligned}
  \end{equation}
  i.e., $\alpha = \rho^2\sign(\rho)$ and $\beta=1/2$.
  
  The limit distribution $H(z)$ in \eqref{eq:HT_assumption} is a centred Gaussian with variance $2\rho^2(1-\rho^2)$, and the penultimate distribution \eqref{eq:penultimate_residual_distribution} is
  \[
    H_x(z)\sim \mathcal N\!\left\{0,\,2\rho^2(1-\rho^2)\left(1+\dfrac{\log x}{2\sqrt{2\rho^2x}}\right)^2\right\}.
  \]
  
  If we write $n^{-1}=\pr(X>u)$, the rate of convergence to the limit distribution is $O\{\log\log n/\sqrt{\log n}\}$ using the ultimate norming in \eqref{eq:HT_convergence}, which is not improved using the penultimate norming in \eqref{eq:HT_convergence_penultimate_norming}. The subasymptotic remainder \eqref{eq:HT_remainder_penultimate} is $O(1/\sqrt{\log n})$.
\end{theorem}

Note that if $b_0(x)=1+\rho x^{1/2}$ and $b_1(x) = 1 + \rho x^{1/2-1/(4x)}$, then the results of Theorem~\ref{thm:gaussian} also hold for all $\rho$, i.e., including $\rho = 0$, with the key change being that the variances of $H$ and $H_x$ no longer have the $\rho^2$ term.

The penultimate norming $a_1(\cdot)$, $b_1(\cdot)$ can be used to assess the goodness-of-fit at a finite level. By replacing $x$ by the threshold $u$ in the log-term of \eqref{eq:gaussian_penultimate_normings}, we derive a second-order approximation for $\alpha=a_1(x)/x$ of the form
\begin{equation}\label{eq:second-order_alpha}
 \alpha_1 = \sign(\rho)\rho^2 + \dfrac{(1-\rho^2)\log u}{2u}.
\end{equation}
Similarly, we derive a second-order approximation for $\beta=\log b_1(x)/\log x$, i.e., 
\begin{equation}\label{eq:second-order_beta}
 \beta_1 = \dfrac{1}{2} - \dfrac{1}{4u}.
\end{equation}

\begin{figure}[htbp]
 \includegraphics[width=\textwidth]{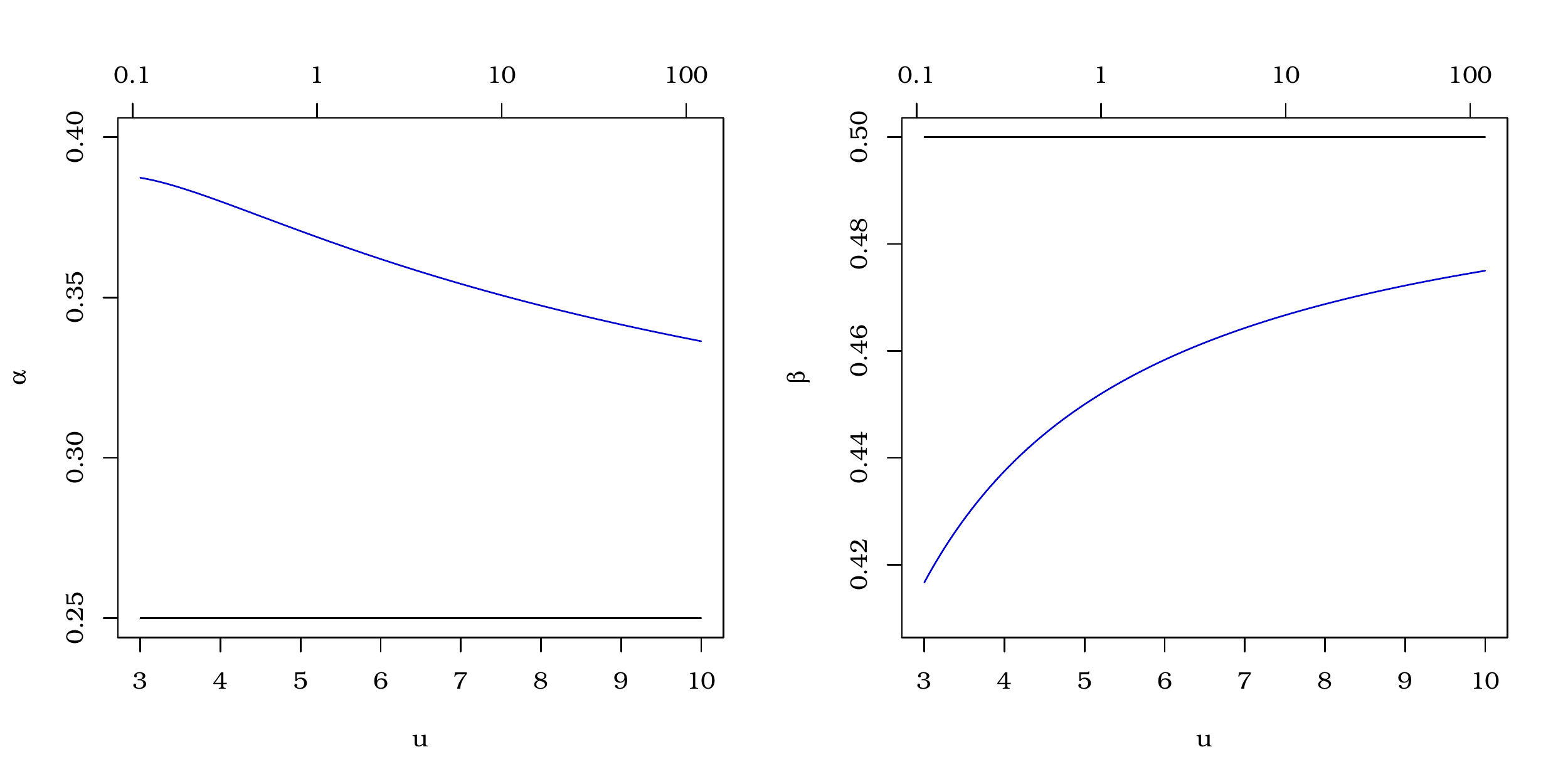}
 \caption{Comparison of first- (black) and second-order (blue) approximations to the Heffernan--Tawn parameters $\alpha$ and $\beta$ for a Gaussian copula with covariance parameter $\rho=0.5$. Lower abscissa on Laplace scale; upper abscissa on the return period scale, in years, assuming daily observations.}\label{fig:gaussian_convergence}
\end{figure}

Figure~\ref{fig:gaussian_convergence} illustrates convergence of the second-order approximations for $\alpha_1$ and $\beta_1$ towards their limits when $\rho=0.5$, and for values of $u$ corresponding to the $97.5\%$ up to the $99.998\%$ Laplace quantile. Convergence is very slow, so  it makes sense to consider second-order approximations when measuring the adequacy of finite-sample estimates. In order to give an idea of the amount of data needed to reach such quantiles, we change the scale of the abscissa to the return period scale, using
\[
  \dfrac{1}{1-F_\text{L}(x)}\times\dfrac{1}{n_Y},
\]
with $F_\text{L}(\cdot)$ the Laplace distribution function, $x$ any quantile on the Laplace scale and $n_Y=365.25$ the number of observations per year.
Figure~\ref{fig:gaussian_convergence} shows that even with the equivalent of more than $100$ years of daily data, the location and scale parameters differ significantly from their asymptotic values.

\subsection{Inverted logistic distribution}
In this section, we consider the bivariate random vector $(X,Y)$ with inverted logistic distribution and Laplace margins \citep{LedfordTawn1997, PapaTawn2016}. Its joint survival distribution function is
\[
 \pr(X>x,Y>y) = 
   \exp\left[-V\left\{\dfrac{-1}{\log\left(\tfrac{1}{2}e^{-x}\right)},\,\dfrac{-1}{\log\left(\tfrac{1}{2}e^{-y}\right)}\right\}\right],\quad x,y>0,
\]
where $V(z,w)=(z^{-1/\gamma}+w^{-1/\gamma})^\gamma$, $0<\gamma\leq 1$, is the exponent measure function of the logistic distribution. Here, $\chi = 0$ and $\bar\chi=2^{1-\gamma}-1$ for $0<\gamma\leq 1$, with $\gamma\to 0$ corresponding to complete dependence and $\gamma = 1$ corresponding to independence.

\begin{theorem}\label{thm:inv_log}
  Let $(X,Y)$ have a bivariate inverted logistic distribution with dependence parameter $0 < \gamma \leq 1$ and Laplace margins. Then the ultimate and penultimate normings \eqref{eq:ultimate_normings} and \eqref{eq:penultimate_normings} for $Y\mid X=x$, with $x$ large, are
  \begin{equation*}
    \begin{aligned}
      a_0(x) &\equiv 0,& b_0(x)&=x^{1-\gamma},\\
      a_1(x) &\equiv-\log 2,& b_1(x)&=x^{1-\gamma},
    \end{aligned}
  \end{equation*}
  so there is no difference in the penultimate form for $b_1(\cdot)$ from $b_0(x)$.
  
  The limit distribution $H(z)$ in \eqref{eq:HT_assumption} is Weibull, specifically $\bar H(z)=\exp(-\gamma z^{1/\gamma})$, and the penultimate distribution $H_x(\cdot)$ in \eqref{eq:penultimate_residual_distribution} is such that
  \begin{equation}\label{eq:invlog_conditional}
    -\log \bar H_x(z)=\gamma z^{1/\gamma} +
         \begin{cases}
            \dfrac{(1-\gamma)(1-\log 2)}{x}z^{1/\gamma}-\dfrac{\gamma(1-\gamma)}{2x}z^{2/\gamma},& 0 < \gamma < 2/3,\\
            \dfrac{1-\log 2}{3x}z^{3/2}-\dfrac{1}{9x}z^3-\dfrac{(\log 2)^2}{8x}\left(4-\dfrac{13\log 2}{3}\right)z^{-3/2},& \gamma=2/3,\\
            -\dfrac{(\log 2)^2}{6\gamma^2}(1-\gamma)\left\{6\gamma+(1-8\gamma)\log 2\right\}x^{3\gamma-3}z^{1/\gamma-3},& 2/3<\gamma< 1.\\
         \end{cases}
  \end{equation}
  When $0<\gamma<2/3$, $H_x(\cdot)$ has finite support
  \[
[0,z_x^H] =     \left[0,\,\left\{\dfrac{x}{1-\gamma}+\dfrac{1-\log 2}{\gamma}\right\}^\gamma\right]
      \longrightarrow\mathbb R_+,\quad x\to\infty,
  \]
  and $1 - H_x(z_x^H)\sim \exp\{-\gamma x/(2-2\gamma)\}$ as $x\to\infty$.
  When $\gamma=2/3$, $H_x(\cdot)$ has finite support, with first-order expansion
  \[
    \left[(12-13\log 2)^{1/3}\left(\dfrac{\log 2}{4}\right)^{2/3} x^{-1/3},\,
    9^{1/3}x^{2/3}\right]\to\mathbb R_+,\quad x\to\infty,
  \]
  and $1-H_x(z_x^H)\sim \exp(-x)$ as $x\to\infty$.
  When $2/3<\gamma<1$, $H_x(\cdot)$ has finite support, with first-order expansion
  \[
    \left[\dfrac{(\log 2)^{2/3}}{\gamma}\left(\dfrac{1-\gamma}{6}\right)^{1/3}\left\{6\gamma+(1-8\gamma)\log 2\right\}^{1/3}x^{\gamma-1},\,
    +\infty\right)\to\mathbb R_+,\quad x\to\infty.
  \]

  If we write $n^{-1}=\pr(X>u)$, the rates of convergence to the limit distribution are $O\{(\log n)^{-1}\}$ using the ultimate norming in \eqref{eq:HT_convergence} and $O\{(\log n)^{\gamma-1}\}$ using the penultimate norming in \eqref{eq:HT_convergence_penultimate_norming}.  The subasymptotic remainder \eqref{eq:HT_remainder_penultimate} behaves like
  \begin{equation}\label{eq:invlog_rates}
    \begin{aligned}
      O\left\{(\log n)^{\alpha-2}\right\},\quad &\alpha\in (0,1/2],&
      O\left\{(\log n)^{3\alpha-3}\right\},\quad &\alpha\in(1/2,2/3),\\
      O\left\{(\log n)^{-4/3}\right\},\quad &\alpha=2/3,&
      O\left\{(\log n)^{-1}\right\},\quad &\alpha\in(2/3,1).
    \end{aligned}
  \end{equation}
\end{theorem}

\begin{figure}
  \centering
  \includegraphics[width=\textwidth]{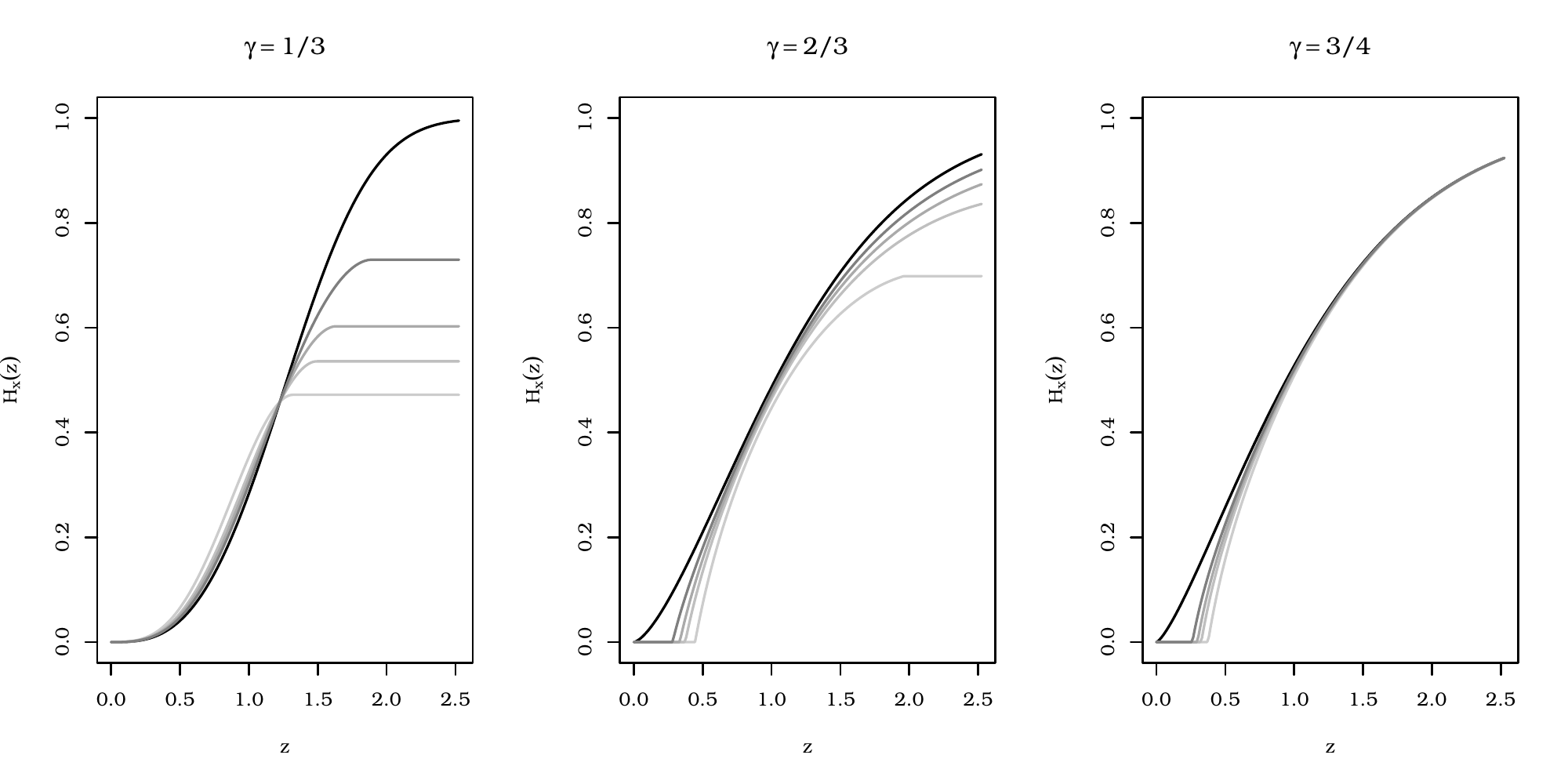}
  \caption{Convergence of the subasymptotic distribution $H_x$ (grey) towards $H$ (black) for the inverted logistic distribution, with $x$ corresponding to the $0.8$, $0.9$, $0.95$ and $0.99$ quantiles (the darker the higher). The panels illustrate the convergence for different values of the dependence parameter, with $\gamma=1/3,2/3,3/4$ from left to right.}
  \label{figure2}
\end{figure}

Figure~\ref{figure2} illustrates the convergence of $H_x$ to $H$ for $\gamma=1/3$, 2/3 and 3/4, and $x$ corresponding to 0.8, 0.9, 0.95 and 0.99 quantiles. It appears that the adequacy of the approximation depends very strongly on $\gamma$.

\subsection{Logistic distribution}
Let $(X,Y)$ have a bivariate logistic distribution with Laplace margins, 
\[
  \pr(X\leq x,\, Y\leq y)=\exp\left[
        -V\left\{\dfrac{-1}{\log\left(1-\tfrac{1}{2}e^{-x}\right)},\,
                 \dfrac{-1}{\log\left(1-\tfrac{1}{2}e^{-y}\right)}\right\}
        \right],\quad x, y>0, 
\]
with
\[
  V(z,w)=\left(z^{-1/\gamma}+w^{-1/\gamma}\right)^\gamma,\quad \gamma\in(0,1].
\]
In the following, we do not consider the case $\gamma = 1$ corresponding to the trivial situation of complete independence. The degree of asymptotic dependence is $\chi = 2-2^\gamma$.
 
\begin{theorem}\label{thm:logistic}
  Let $(X,Y)$ have a bivariate inverted logistic distribution with dependence parameter $0<\gamma\leq 1$ and Laplace margins. Then, the ultimate normings \eqref{eq:ultimate_normings} for $Y\mid X=x$, with $x$ large, are
  \begin{equation*}
    a_0(x)=x,\quad b_0(x)=1,
  \end{equation*}
  and penultimate normings \eqref{eq:penultimate_normings} are identical to $a_0$, $b_0$.
\end{theorem}
 
\section{Penultimate model form}\label{sec:discussion}
The three examples of copula studied in Section~\ref{sec:examples} all have very different extremal dependence features, yet all are in the following, rather general, class of penultimate forms for the norming functions
\[
a_1(x)= \left(\alpha+\dfrac{{\cal L}_{a}(x)}{x^{\gamma_{a}}}\right) x,\qquad
\log b_1(x)= \left(\beta+\dfrac{{\cal L}_b(x)}{x^{\gamma_{b}}}\right) \log x 
\]
where $\alpha$ and $\beta$ are from of the first order norming functions $a_0$ and $b_0$ respectively, $\gamma_{a}>0$, $\gamma_{b}\ge 0$ and ${\cal L}_{a}$, ${\cal L}_{b}$ are functions that are slowly varying at $\infty$. For statistical modelling purposes we need to be more precise about the slowly varying function, and as in past practice, e.g., \citet{LedfordTawn1996}, we fix these functions above some threshold $u$ to be constants, i.e.,
\[
{\cal L}_{a}(x)=\delta_a, \qquad {\cal L}_{b}(x)=\delta_b,\quad x>u.
\]
In practice this may still be rather over-parameterised, given how difficult second-order effects are to estimate, and so in practice it may be sufficient to fix $\gamma_a=\gamma_b=1$. Choices like this have been used in penultimate modelling in univariate cases \citep{HallWelsh1985}.

Next consider the choice of the limit distribution $H$. Although our examples have shown we can get improved penultimate forms for $H$, using $H_x$, these improvements do not always lead to faster convergence to the limit distribution for certain ranges of parameters of the underlying model. Furthermore, as the  current approach is to use a non-parametric approach to estimate $H$, it is suggested that this is continued and $H$ is assumed independent of $x$.

In conclusion, based on a series of examples, we have suggested a class of penultimate models for the \cite{HeffernanTawn2004} conditional extremes model. This new class of models improves convergence rates relative to the limit model, and thus is likely to lead to better statistical inference. The extension is parsimonious with, in its most simple form, the introduction of two additional parameters.

\section*{Acknowledgements}
This research was partially supported by the Swiss National Science Foundation.

\appendix
\section{Appendix}
\subsection{Proof of Theorem~\ref{thm:gaussian}}
The details of the proof for the asymptotic quantities $a_0(\cdot)$, $b_0(\cdot)$ and $H(z)$ can be found in \citet{HeffernanTawn2004}. We follow a similar path to derive the penultimate approximations and use Mill's ratio to get a tail approximation for $v$ and $x$ on the normal and Laplace scales, respectively. Since the Laplace distribution is symmetric, we focus on the right tail:
\begin{align}\label{eq:Gaussian_second-order_y}
 x &\sim -\log\left\{2\dfrac{\varphi(v)}{v}\right\}
      =-\log 2+\log v+\dfrac{1}{2}\log(2\pi)+\dfrac{1}{2}v^2,\\
 v &\sim \sqrt{2x},\nonumber
\end{align}
for large $x$ and $v$.
In order to get a second-order approximation for $v$, we define a small quantity $\varepsilon>0$ and set $v\sim \sqrt{2x}(1+\varepsilon)$; plugging this  back in \eqref{eq:Gaussian_second-order_y}, we get
\begin{align}\nonumber
 x &\sim -\log 2 + \log\left\{\sqrt{2x}(1+\varepsilon)\right\}
           + \dfrac{1}{2}\log(2\pi)+\dfrac{1}{2}\left\{\sqrt{2x}(1+\varepsilon)\right\}^2\\
     &\sim -\log 2 +\dfrac{1}{2}\log\left(\sqrt{2x}\right)+\dfrac{1}{2}\log(2\pi) +x +2x\varepsilon,
\end{align}
and this yields
\begin{equation}
     \label{eq:v_second_order}
     v= \sqrt{2x}+\dfrac{2\log 2-\log(2x)-\log(2\pi)}{2\sqrt{2x}} + O\left\{\dfrac{(\log x)^2}{x^{3/2}}\right\}.
\end{equation}

We want the conditional distribution of $W\mid V=v$ to be well-behaved in its upper tail when $v\rightarrow \infty$. We have $\pr(W-\rho V\leq z\mid V=v)=\Phi(z/\sqrt{1-\rho^2})$. On the Laplace scale, $W$ is transformed to
\begin{equation}\label{eq:HT_model_penultimate}
 Y = a_1(x)+b_1(x)Z,
\end{equation}
with $a_1(x)$ and $b_1(x)>0$ norming functions to be determined, and $Z$ a random variable with a fixed distribution non-degenerate at $+\infty$.
We derive these penultimate norming functions by writing $a_1(x)=a_0(x)(1+\varepsilon)=\rho^2x(1+\varepsilon)$ in  \eqref{eq:HT_model_penultimate}, and using \eqref{eq:v_second_order}. We get the second-order approximation by first expanding
\begin{equation}\label{eq:Gaussian_first_expansion}
\begin{aligned}
 W \sim& \sqrt{2\rho^2x(1+\varepsilon)+2b(x)Z}-\dfrac{\log\pi+\log\{\rho^2x(1+\varepsilon)+b(x)Z\}}{2\sqrt{2\rho^2x(1+\varepsilon)+2b(x)Z}}\\
     \sim& \rho\sqrt{2x(1+\varepsilon)}\left\{1+\dfrac{b(x)Z}{2\rho^2x(1+\varepsilon)}\right\}\\
           &-\left[\log\pi+\log\left\{\rho^2x(1+\varepsilon)\right\}+\dfrac{b(x)Z}{\rho^2x(1+\varepsilon)}\right]\dfrac{1}{2\sqrt{2\rho^2x(1+\varepsilon)}}\left\{1-\dfrac{b(x)Z}{2\rho^2x(1+\varepsilon)}\right\}.
\end{aligned}
\end{equation}
We can use this expression for $W$ in $W-\rho V\mid V=v$, in which we keep only the terms that are not functions of $Z$, as they are disconnected from $a_1(\cdot)$, and we get
\begin{equation*}
  \rho\sqrt{2x(1+\varepsilon)}-\dfrac{\log\left\{\rho^2x(1+\varepsilon)\right\}}{2\sqrt{2\rho^2x(1+\varepsilon)}}
  - \rho\sqrt{2x} +\rho\dfrac{\log x}{2\sqrt{2x}}+O\left(x^{-1/2}\right).
\end{equation*}
Expanding further, we arrive at
\[
   \rho\sqrt{2x}\left(1+\dfrac{\varepsilon}{2}\right) - \dfrac{\log(\rho^2x)+\varepsilon}{2\sqrt{2\rho^2x}}\left(1-\dfrac{\varepsilon}{2}\right)
   - \rho\sqrt{2x} + \rho\dfrac{\log x}{2\sqrt{2x}}+O\left(x^{-1/2}\right),
\]
and cancellation of the leading term yields
\begin{equation}\label{eq:second-order_eps}
   \varepsilon= \dfrac{(1-\rho^2)\log x}{2\rho^2x},
\end{equation}
or equivalently $a_1(x)=\rho^2x+\tfrac{(1-\rho^2)}{2}\log x$.

The penultimate scale function $b_1(\cdot)$ stems from the $Z$-terms in \eqref{eq:Gaussian_first_expansion}, namely
\begin{equation}\label{eq:Gaussian_second-order_beta}
   \dfrac{b(x)Z}{\rho\sqrt{2x(1+\varepsilon)}}+\dfrac{\log\left\{\rho^2x(1+\varepsilon)\right\}b(x)Z}{2\left\{2\rho^2x(1+\varepsilon)\right\}^{3/2}} - \dfrac{b(x)Z}{\left\{2\rho^2x(1+\varepsilon)\right\}^{3/2}},
\end{equation}
which we expand as
\begin{align*}
 b(x) &\sim \rho\sqrt{2x(1+\varepsilon)}\left[1+\dfrac{\log\left\{x(1+\varepsilon)\right\}}{4\rho^2x(1+\varepsilon)} - \dfrac{1}{2\rho^2x(1+\varepsilon)}\right]^{-1}\\
     & \sim    \rho\sqrt{2x}\left(1+\dfrac{\varepsilon}{2}\right)\left\{1-\dfrac{\log x}{4\rho^2x}(1-\varepsilon)\right\}\\
     & \sim    \rho\sqrt{2x}-\dfrac{\log x}{2\rho\sqrt{2x}}+\rho\sqrt{2x}\dfrac{\varepsilon}{2}.
\end{align*}
Substituting $\varepsilon$ into \eqref{eq:second-order_eps} gives
\begin{equation*}
 b(x)  \sim \rho\sqrt{2x}-\dfrac{\log x}{2\rho\sqrt{2x}}+\dfrac{(1-\rho^2)\log x}{2\rho\sqrt{2x}}
         \sim x^{1/2-1/(4x)}=b_1(x).
\end{equation*}

We now compute the penultimate distribution $H_x(z)$ by substituting the expression for $\varepsilon$ in \eqref{eq:second-order_eps}, writing $A(x)+B(x)Z\sim H_x(z)$, with
\[
  \rho\sqrt{2x}\dfrac{(1-\rho^2)\log x}{4\rho^2x}
  -\dfrac{\log(\rho^2x)}{2\sqrt{2\rho^2x}}
  -\dfrac{(1-\rho^2)\log x}{4\rho^2x\sqrt{2\rho^2x}}
  +\dfrac{(1-\rho^2)\log(\rho^2x)\log x}{4\sqrt{2\rho^2x}(2\rho^2x)}
  +\rho\dfrac{\log x}{2\sqrt{2x}},
\]
which equals
\begin{equation}\label{eq:non-z_part}
     -\dfrac{\log(\rho)}{\sqrt{2\rho^2x}}
  -\dfrac{(1-\rho^2)\log(x)}{2(2\rho^2x)^{3/2}}
  +\dfrac{(1-\rho^2)(\log x)^2}{4(2\rho^2x)^{3/2}}
  +\dfrac{(1-\rho^2)\log(\rho)}{2(2\rho^2x)^{3/2}}
  \sim -\dfrac{\log(\rho)}{\sqrt{2\rho^2x}} = A(x),
\end{equation}
and
\[
  \dfrac{b(x)Z}{\sqrt{2\rho^2x}}
  \left\{(1-\varepsilon/2)
      +\dfrac{\log(\rho^2x)+\varepsilon}{2\sqrt{2\rho^2x}}(1-\varepsilon)
      -\dfrac{1-3\varepsilon/2}{2\rho^2x}\right\},
\]
which equals
\begin{multline}\label{eq:z-part}
  \dfrac{b(x)Z}{\sqrt{2\rho^2x}}
         \left\{1-\dfrac{\varepsilon}{2}
      +\dfrac{\log x}{2\sqrt{2\rho^2 x}}(1-\varepsilon)
      +\dfrac{\varepsilon}{2\sqrt{2\rho^2 x}}
      +\dfrac{\log \rho}{\sqrt{2\rho^2 x}}(1-\varepsilon)
      -\dfrac{1}{2\rho^2 x} + O(\varepsilon^2)\right\}\\
  =     Z+ Z\dfrac{\log x}{2\sqrt{2\rho^2x}}
      +Z\dfrac{\log \rho}{\sqrt{2\rho^2x}}
      +Z\times O\left(\dfrac{\log x}{x}\right).
\end{multline}
Taking the leading and penultimate terms in \eqref{eq:non-z_part} and \eqref{eq:z-part}, which corresponds to setting $B(x)=1 + (8\rho^2 x)^{-1/2}\log x$, we find that $H_x(z)=H_x^{(1)}(z)$ is a centred normal distribution with variance
\[
   2\rho^2(1-\rho^2)\left(1+\dfrac{\log x}{2\sqrt{2\rho^2x}}\right)^2,
\]
i.e., it has larger variance than the asymptotic $H(z)$ found in \citet{HeffernanTawn2004}. If we consider the antepenultimate terms in \eqref{eq:non-z_part} and \eqref{eq:z-part}, we get $H_x(z)=H^{(2)}_x(z)$ as
\[
  \mathcal N\left\{-\dfrac{\log(\rho)}{\sqrt{2\rho^2x}},\,
   2\rho^2(1-\rho^2)\left(1+\dfrac{\log x+2\log \rho}{2\sqrt{2\rho^2x}}\right)^2\right\}.
\]
We now give the order of convergence of \eqref{eq:HT_remainder_penultimate}
with the penultimate approximations $a_1(\cdot)$ and $b_1(\cdot)$ that we derived. After a marginal transform in order to get $n^{-1}=\pr(X>u)$, we get the rate of convergence for $H_x^{(1)}(z)$, namely $O(1/\sqrt{\log n})$, which improves on the version of $H(z)$ with first-order approximation for the normalising functions, for which the rate of convergence is of order $\log \log n/\sqrt{\log n}$. Taking $H_x^{(2)}(z)$ improves even more on $H(z)$, as its rate of convergence in \eqref{eq:HT_remainder_penultimate} is of order $O(\log\log n/\log n)$.

\subsection{Proof of Theorem~\ref{thm:inv_log}}
We start by computing the conditional survival distribution of $Y\mid X=x$ for large $x$ and deriving a tail approximation to it. We have, for non-negative $x$ and $y$,
\begin{align*}
  \pr(Y > y&\mid X =x)\\
  &= 2e^{x}\times\dfrac{\partial \pr(X>x,\,Y>y)}{\partial x}\\
  &= -2\exp\left\{x-V\left(\dfrac{1}{x+\log 2},\,\dfrac{1}{y+\log 2}\right)\right\}
    V_1\left(\dfrac{1}{x+\log 2},\,\dfrac{1}{y+\log 2}\right)
    (x+\log 2)^{-2},
\end{align*}
where the partial derivative of the exponent measure is
\[
  V_1(x,y) = -\left(x^{-1/\gamma}+y^{-1/\gamma}\right)^{\gamma-1} x^{-1/\gamma-1}.
\]

To ease  the following developments, we examine the log-survival conditional probability
\begin{align*}
  \log\pr\left(Y > y\mid\right.&\left. X =x\right)\\
  =&\log 2+x-\left[\left\{x(1+x^{-1}\log 2)\right\}^{1/\gamma}
                    + \left\{y(1+y^{-1}\log2)\right\}^{1/\gamma}\right]^\gamma\\
    &+ (\gamma-1)\log\left[\left\{x(1+x^{-1}\log 2)\right\}^{1/\gamma}+\left\{y(1+y^{-1}\log 2)\right\}^{1/\gamma}\right]\\
    &+\dfrac{1-\gamma}{\gamma}\log\left\{x(1+x^{-1}\log 2)\right\}\\
  \approx & \log 2 + x - x\left[
        1+\dfrac{\log 2}{\gamma x}
        +\dfrac{1-\gamma}{2\gamma^2}\dfrac{(\log 2)^2}{x^2}
      +\left(\dfrac{y}{x}\right)^{1/\gamma}\left\{1+\dfrac{\log 2}{\gamma y}
        +\dfrac{1-\gamma}{2\gamma^2}\dfrac{(\log 2)^2}{y^2}\right\}
      \right]^\gamma\\
      &+(\gamma-1)\left[\dfrac{1-\gamma}{2\gamma^2}\dfrac{(\log 2)^2}{x^2}
        +\left(\dfrac{y}{x}\right)^{1/\gamma}\left\{
            1+\dfrac{\log 2}{\gamma y}+\dfrac{1-\gamma}{2\gamma^2}\dfrac{(\log 2)^2}{y^2}
      \right\}\right],
\end{align*}
for large $x$ and $y$. We can expand further, using the fact that $x$ and $y$ are positively associated and asymptotically independent, so large values of $x$ occur with large values of $y$ with large ratio $x/y$, so the log conditional probability can be approximated as follows,
\begin{multline*}
  -x\Bigg[\dfrac{1-\gamma}{2\gamma}\dfrac{(\log 2)^2}{x^2}
      +\left(\dfrac{y}{x}\right)^{1/\gamma}\left\{
        \gamma + \dfrac{\log 2}{y}+\dfrac{1-\gamma}{2\gamma}\dfrac{(\log 2)^2}{y^2}\right\}\\
      +\dfrac{\gamma(\gamma-1)}{2}\left\{\dfrac{\log 2}{\gamma x}
      +\left(\dfrac{y}{x}\right)^{1/\gamma}\left(1+\dfrac{\log 2}{\gamma y}\right)
      \right\}^2\Bigg]\\
      -\dfrac{(1-\gamma)^2}{2\gamma}\dfrac{(\log 2)^2}{x^2}
      -(1-\gamma)\left(\dfrac{y}{x}\right)^{1/\gamma}\left\{
      1+\dfrac{\log 2}{\gamma y}+\dfrac{1-\gamma}{2\gamma^2}\dfrac{(\log 2)^2}{y^2}
      \right\}.
\end{multline*}

For $Y=a(x)+b(x)Z$, the first order behaviour is cancelled by choosing $a_0(x)\equiv 0$ and $b_0(x)=x^{1-\gamma}$ \citep{HeffernanTawn2004}. We find $a_1(\cdot)$ by setting $a_1(x)=\varepsilon$, with $\varepsilon=\varepsilon(x)=o(x^{1-\gamma})$, namely with $Y=\varepsilon+x^{1-\gamma}Z$,
\begin{align*}
  \nonumber
  \log&\pr(Y > y\mid X =x)\\ \nonumber
  \approx & -\dfrac{1-\gamma}{2\gamma}(\log 2)^2x^{-1}
  -x^{1-1/\gamma}\left(\varepsilon+x^{1-\gamma}z\right)^{1/\gamma}\\
  &\hspace{20ex}\times\left\{
      \gamma+(\log 2)\left(\varepsilon+x^{1-\gamma}z\right)^{-1}+\dfrac{1-\gamma}{2\gamma}(\log 2)^2\left(\varepsilon+x^{1-\gamma}z\right)^{-2}
      \right\}\\ \nonumber
  &-\dfrac{(\gamma-1)(\log 2)^2}{2\gamma x}
  -\dfrac{\gamma(\gamma-1)}{2}x^{1-2/\gamma}\left(\varepsilon+x^{1-\gamma}z\right)^{1/\gamma}\left\{
      1+\dfrac{\log 2}{\gamma}\left(\varepsilon+x^{1-\gamma}z\right)^{-1}\right\}\\ \nonumber
  &+(1-\gamma)(\log 2)x^{-1/\gamma}\left(\varepsilon+x^{1-\gamma}z\right)^{1/\gamma}\left\{
      1+\dfrac{\log 2}{\gamma}\left(\varepsilon+x^{1-\gamma}z\right)^{-1}\right\}
  -\dfrac{(1-\gamma)^2}{2\gamma}(\log 2)^2x^{-2}\\ \nonumber
  &-(1-\gamma)x^{-1/\gamma}\left(\varepsilon+x^{1-\gamma}z\right)^{1/\gamma}\left\{
  1+\dfrac{\log 2}{\gamma}\left(\varepsilon+x^{1-\gamma}z\right)^{-1}
      +\dfrac{1-\gamma}{2\gamma^2}(\log 2)^2\left(\varepsilon+x^{1-\gamma}z\right)^{-2}\right\}\\
  = & -z^{1/\gamma}\left\{1+\dfrac{\varepsilon}{\gamma}x^{\gamma-1}z^{-1}
        +\dfrac{1-\gamma}{2\gamma^2}\varepsilon^2x^{2\gamma-2}z^{-2}
        +\dfrac{(1-\gamma)(1-2\gamma)}{6\gamma^3}\varepsilon^3x^{3\gamma-3}z^{-3}
        +O\left(x^{4\gamma-4}\right)\right\}\\
        &\times\Bigg\{\gamma+(\log 2)\left(x^{\gamma-1}z^{-1}
        -\varepsilon x^{2\gamma-2}z^{-2}+\varepsilon^2 x^{3\gamma-3}z^{-3}\right)\\
        &\hspace{20ex}+\dfrac{1-\gamma}{2\gamma}(\log 2)^2\left(x^{2\gamma-2}z^{-2}
        -2\varepsilon x^{3\gamma-3}z^{-3}\right)+O\left(x^{4\gamma-4}\right)
        \Bigg\}\\
        &+\dfrac{\gamma(1-\gamma)}{2}z^{2/\gamma}x^{-1}
        +(1-\gamma)(\log 2)x^{-1}z^{1/\gamma}
        -(1-\gamma)x^{-1}z^{1/\gamma}+O\left(x^{\gamma-2}\right).
\end{align*}
Expanding this expression and rearranging the terms yields
\begin{equation}\label{eq:invlog_expansion}
\begin{aligned}
  &-\gamma z^{1/\gamma}-(\log 2+\varepsilon)x^{\gamma-1}z^{1/\gamma-1}
  +\left\{
      (1-\gamma)(\log 2-1)z^{1/\gamma}+\dfrac{\gamma(1-\gamma)}{2}z^{2/\gamma}
  \right\}x^{-1}\\
  &+\left\{
      (\log 2)\varepsilon z^{1/\gamma-2}-\dfrac{\varepsilon\log 2}{\gamma}z^{1/\gamma-2}
      -\dfrac{1-\gamma}{2\gamma}\varepsilon^2 z^{1/\gamma-2}
      -\dfrac{1-\gamma}{2\gamma}(\log 2)^2 z^{1/\gamma-2}
  \right\}x^{2\gamma-2}\\
  &-\left\{
      \varepsilon^2 -\dfrac{1-\gamma}{\gamma}(\log 2)^2\varepsilon
      -\dfrac{\varepsilon^2}{\gamma}+\dfrac{1-\gamma}{2\gamma^2}(\log 2)^2\varepsilon
      +\dfrac{1-\gamma}{2\gamma^2}(\log 2)\varepsilon^2
      +\dfrac{(1-\gamma)(1-2\gamma)}{6\gamma^2}\varepsilon^3
  \right\}z^{1/\gamma-3}x^{3\gamma-3}\\
  &+O\left(x^{\max\{\gamma-2,4\gamma-4\}}\right).
\end{aligned}
\end{equation}
We obtain $\varepsilon=-\log 2$, or equivalently $a_1(x)\equiv -\log 2$, cancelling the leading term in $x$ in \eqref{eq:invlog_expansion}. 
Higher-order terms imply different powers of $Z$, hence further approximation of the norming functions is infeasible because of the linearity in $Z$ stemming from the location-scale norming.

Plugging $a_1(x)$ into \eqref{eq:invlog_expansion} yields 
\begin{multline*}
  -\gamma z^{1/\gamma}
  + (1-\gamma)\left\{
  \log 2 - 1+\dfrac{\gamma}{2}z^{1/\gamma}\right\}z^{1/\gamma}x^{-1}\\
  +\dfrac{(1-\gamma)(\log 2)^2}{6\gamma^2}\left\{6\gamma+(1-8\gamma)\log 2\right\}z^{1/\gamma-3}x^{3\gamma-3}
  +O\left(x^{\max\{\gamma-2,4\gamma-4\}}\right).
\end{multline*}
The various expressions for $H_x(\cdot)$ depending on the value of $\gamma$ directly follow.

We now derive the support of $H_x(\cdot)$ when $0<\gamma<2/3$. A necessary condition for $H_x(\cdot)$ to be well-defined is that the density $h_x(\cdot)=H_x'(\cdot)$ is non-negative, that is
\begin{align*}
  &\bar H_x(z)\left\{z^{1/\gamma-1}+\dfrac{(1-\gamma)(1-\log 2)}{\gamma x}z^{1/\gamma-1}-\dfrac{1-\gamma}{x}z^{2/\gamma-1}\right\} \geq 0
  \end{align*}
yielding 
\begin{align*}
& \dfrac{1-\gamma}{x}z^{2/\gamma} \leq z^{1/\gamma}\left\{1+\dfrac{(1-\gamma)(1-\log 2)}{\gamma x}\right\}
\end{align*}
and 
\begin{align*}
& z\leq \left(\dfrac{x}{1-\gamma}+\dfrac{1-\log 2}{\gamma}\right)^\gamma,
\end{align*}
and the upper bound is also the upper endpoint $z^H$; the value of $H_x(z^H)$ follows directly. The lower endpoint is attained when the exponent in \eqref{eq:invlog_conditional} vanishes, in other terms when
\[
  z^{1/\gamma}\left(2\dfrac{1-\log 2}{\gamma}+\dfrac{2x}{1-\gamma}
  -z^{1/\gamma}\right)=0,
\]
for which the root of interest is $z=0$, which concludes this part of the proof.

The case $\gamma=2/3$ is treated similarly, as we require the derivative of $H_x(z)$ to be non-negative,
\begin{align}
            & \dfrac{1}{3x}z^2-z^{1/2}-\dfrac{3}{16x}(\log 2)^2\left(4-\dfrac{13}{3}\log 2\right)z^{-5/2}\leq 0, \label{eq:upper_endpoint}
\end{align}
or equivalently
\begin{align}\nonumber
& w^3-(3x+a)w^2-c\leq 0,& z>0,x>0,
\end{align}
with $w=z^{3/2}$, $a=3(1-\log 2)/2$ and $c=3(\log 2)^2(12-13\log 2)/16$. When $x\to\infty$, we have $w\to\infty$, and a leading term is given by $  w-(3x+a)\leq 0$, implying that $w\leq 3x+a$. 
In order to ensure that we are not missing a term in this approximation, 
we consider \eqref{eq:upper_endpoint} with $w= 3x+a+\delta$,  $0<\delta=\delta(x)=O(x)$, as follows,
\begin{equation*}
  (3x+a+\delta)^3-(3x+a)(3x+a+\delta)^2-c
  =\delta(3x+a)^2+2\delta^2(3x+a)+\delta^3-c\leq 0,\quad x>0.
\end{equation*}
For this inequality to hold, we need at least $9x^2\delta\leq c$, $x>0$, i.e., $\delta=O(x^{-2})$. We conclude that $z^H=\{3x+3(1-\log 2)/2\}^{2/3}$ is an approximation of the upper endpoint when $\gamma=2/3$, with
\[
  \bar H_x\left(z^H\right) = \exp\left\{-\dfrac{1}{x}\left(x+\dfrac{1-\log 2}{2}\right)^2+O\left(x^{-2}\right)\right\}\sim \exp(-x)
  \to 0,\quad x\to \infty.
\]

The lower endpoint is computed by finding an approximation to the root of interest of the exponent in \eqref{eq:invlog_conditional}, or equivalently with $w=z^{3/2}$,
\[
  w^3-\left\{3(1-\log 2)+6x\right\}w^2+\dfrac{3}{8}(\log 2)^2(12-13\log 2)=0,
\]
for which we know $w\to 0$ when $x\to \infty$, leading to the approximation
\begin{equation}\label{eq:lower_endpoint}
  (3x+a)w^2+c=0\\
  \implies w=\left(\dfrac{c}{3x+a}\right)^{1/2}
      \approx \sqrt{c}\left\{\dfrac{1}{\sqrt{3x}}-\dfrac{a}{2(3x)^{3/2}}\right\},\quad w>0,x>0.
\end{equation}
Consider $w=\sqrt{c/(6x+a)}+\delta$, $0<\delta=\delta(x)=O(1/\sqrt{x})$, in order to confirm that \eqref{eq:lower_endpoint} is a sensible approximation as follows,
\[
  \left(\dfrac{c}{3x+a}\right)^{3/2}+3\dfrac{c\delta}{3x+a}+3\left(\dfrac{c}{3x+a}\right)^{1/2}\delta^2+\delta^3-(3x+a)\left\{2\delta\left(\dfrac{c}{3x+a}\right)^{1/2}+\delta^2\right\}=0,
\]
and expanding the expressions in brackets yields
\begin{multline*}
  c^{3/2}\left(3x\right)^{-3/2}+c\delta x^{-1}
  +\sqrt{3c}\delta^2 x^{-1/2}\left(1-\dfrac{a}{6x}\right)
  +\delta^3\\
  -2\sqrt{3c}\delta x^{1/2}\left(1+\dfrac{a}{6x}-\dfrac{a^2}{72x^2}\right)
  -\delta^2a-6\delta^2 x+O\left(x^{-5/2}\right)=0.
\end{multline*}
From this we observe that we require $\delta=o(x^{-1/2})$ for the equality to hold as $x\to\infty$, so we can simplify further and get
\[
  \delta^2\left(3x+a\right)+\delta\left(2\sqrt{3c}x^{1/2}
  +\dfrac{a\sqrt{3c}}{3}x^{-1/2}\right)-c^{3/2}\left(3x\right)^{-3/2} = 0,
\]
which we solve in $\delta$. We get an approximate square root discriminant
\[
  x^{1/2}\left(2\sqrt{3c}+\dfrac{a\sqrt{3c}}{3}x^{-1}+\dfrac{c}{3}x^{-3/2}\right),
\]
from which we compute the approximate root of interest
\[
  \delta = \dfrac{c}{6}x^{-1}\left(6x+2a\right)^{-1}\sim \dfrac{c}{36}x^{-2}.
\]
This ends the proof for the support of $H_x(z)$ when $\gamma=2/3$.

In the case when $\gamma\in(2/3,1)$, we require
\[
  \dfrac{\partial\bar H_x(z)}{\partial z}
  = \bar H_x(z)\left[-z^{1/\gamma-1}+\left(\dfrac{1}{\gamma}-3\right)\dfrac{(\log 2)^2}{6\gamma^2}(1-\gamma)\left\{6\gamma+(1-8\gamma)\log 2\right\}x^{3\gamma-3}z^{1/\gamma-4}\right]\leq 0,
\]
which is true for all $z>0$ and $x>0$. The density is well-defined and we can verify that the upper endpoint of $H_x(z)$ is $+\infty$. We work out the lower endpoint by considering the exponent in \eqref{eq:invlog_conditional}, with
\begin{multline*}
  \gamma z^{1/\gamma} -\dfrac{(\log 2)^2}{6\gamma^2}(1-\gamma)\left\{6\gamma+(1-8\gamma)\log 2\right\}x^{3\gamma-3}z^{1/\gamma-3}\\
  =\gamma z^{1/\gamma}\left[1-\dfrac{(\log 2)^2}{6\gamma^3}(1-\gamma)\left\{
  6\gamma+(1-8\gamma)\log 2\right\}x^{3\gamma-3}z^{-3}\right],
\end{multline*}
which vanishes when $z=0$, and 
\begin{equation}\label{eq:invlog_lower_endpoint_abig}
  z=\dfrac{(\log 2)^{2/3}}{\gamma}\left[\dfrac{1-\gamma}{6}\left\{6\gamma+(1-8\gamma)\log 2\right\}\right]^{1/3}x^{\gamma-1},\quad x>0.
\end{equation}
The root of interest is \eqref{eq:invlog_lower_endpoint_abig}, giving the desired result.

We now give the convergence rate of \eqref{eq:HT_remainder_penultimate} using the penultimate approximation $a_1(\cdot)$. The convergence rate is linked with the value of the dependence parameter $\gamma$ and can be found from \eqref{eq:invlog_expansion}.
The powers of $x$ of interest appearing in \eqref{eq:invlog_expansion} are $-1$, $3\gamma-3$, $\gamma-2$, $4\gamma-4$, depending on the precise value of $\gamma$.
For $\gamma \in (0,1/2]$, convergence is the fastest, as subtraction of $H_x(z)$ in this case removes terms in $x^{-1}$, so we conclude that \eqref{eq:HT_remainder_penultimate} has a leading term in $x^{\gamma-2}$. Similarly we conclude that the order of convergence for $\gamma\in(1/2,2/3)$ is $x^{3\gamma-3}$. For $\gamma=2/3$, the $\gamma-2$ and $4\gamma-4$ powers coincide and give $x^{-4/3}$ as the leading term. When $\gamma\in(2/3,1)$, convergence is slowest with $x^{-1}$ as the leading term.

\subsection{Proof of Theorem~\ref{thm:logistic}}
We now focus on the conditional probability
\begin{equation}\label{eq:log_conditional_start}
  \begin{aligned}
  \pr(Y\leq y\mid X=x)
  =&\; 2\exp\left[
        x
        -V\left\{\dfrac{-1}{\log\left(1-\tfrac{1}{2}e^{-x}\right)},\,
                  \dfrac{-1}{\log\left(1-\tfrac{1}{2}e^{-y}\right)}\right\}
      \right]\\
      &\;\times\left[
          -V_1\left\{\dfrac{-1}{\log\left(1-\tfrac{1}{2}e^{-x}\right)},\,
                      \dfrac{-1}{\log\left(1-\tfrac{1}{2}e^{-y}\right)}\right\}
      \right]
      \dfrac{d}{dx}\left\{-\dfrac{1}{\log\left(1-\tfrac{1}{2}e^{-x}\right)}\right\}.
    \end{aligned}
\end{equation}
We can approximate the last term in \eqref{eq:log_conditional_start}, for large $x$,
\begin{equation}\label{eq:log_derivee_interne}
\begin{aligned}
  \dfrac{d}{dx}\left\{-\dfrac{1}{\log\left(1-\tfrac{1}{2}e^{-x}\right)}\right\}
  &=       \dfrac{1}{\left\{\log\left(1-\tfrac{1}{2}e^{-x}\right)\right\}^2}
            \dfrac{1}{\left(1-\tfrac{1}{2}e^{-x}\right)}
            \dfrac{1}{2}e^{-x}\\
  &\approx \left(\dfrac{1}{2}e^{-x}+\dfrac{1}{8}e^{-2x}\right)^{-2}
            \left(\dfrac{1}{2}e^{-x}+\dfrac{1}{4}e^{-2x}\right)\\
  &\approx 2e^{x}-\dfrac{1}{8}e^{-x}.
\end{aligned}
\end{equation}
The partial derivative of $V(\cdot,\cdot)$ in \eqref{eq:log_conditional_start} can be approximated as
\begin{equation}\label{eq:log_spectral_partial_derivative}
\begin{aligned}
  V_1&\left\{-\dfrac{1}{\log\left(1-\tfrac{1}{2}e^{-\gras y}\right)}\right\}\\
    \approx&-\left\{\left(\dfrac{1}{2}e^{-x}+\dfrac{1}{8}e^{-2x}\right)^{1/\gamma}
            +\left(\dfrac{1}{2}e^{-y}+\dfrac{1}{8}e^{-2y}\right)^{1/\gamma}
      \right\}^{\gamma-1}
      \left(\dfrac{1}{2}e^{-x}+\dfrac{1}{8}e^{-2y}\right)^{1/\gamma+1}\\
    \approx& -\dfrac{1}{4}\left\{
        e^{-x/\gamma}+\dfrac{1}{4\gamma}e^{-x(1+1/\gamma)}
        +\dfrac{1-\gamma}{32\gamma^2}e^{-x(2+1/\gamma)}
        +e^{-y/\gamma}+\dfrac{1}{4\gamma}e^{-y(1+1/\gamma)}
        +\dfrac{1-\gamma}{32\gamma^2}e^{-y(2+1/\gamma)}
        \right\}^{\gamma-1}\\
      &\times
          e^{-x(1+1/\gamma)}\left(1+\dfrac{1}{4}e^{-x}\right)^{1+1/\gamma},
\end{aligned}
\end{equation}
for large $x$. From  \eqref{eq:log_derivee_interne} and \eqref{eq:log_spectral_partial_derivative}, the log conditional probability \eqref{eq:log_conditional_start} is
\begin{equation}\label{eq:log_conditional_cont}
\begin{aligned}
  \log\{\pr(Y\leq y&\mid X=x)\}\\
  =&\; \log 2+x-\left\{
              \left(\dfrac{1}{2}e^{-x}\right)^{1/\gamma}
              \left(1+\dfrac{1}{4\gamma}e^{-x}\right)
            +\left(\dfrac{1}{2}e^{-y}\right)^{1/\gamma}
              \left(1+\dfrac{1}{4\gamma}e^{-y}\right)
            \right\}^\gamma\\
      &\;-\log 4-(1-\gamma)\log\left\{
              e^{-x/\gamma}+\dfrac{1}{4\gamma}e^{-x(1+1/\gamma)}
              +e^{-y/\gamma}+\dfrac{1}{4\gamma}e^{-y(1+1/\gamma)}
                    \right\}\\
      &\;-\dfrac{1+\gamma}{\gamma}x+\dfrac{1+\gamma}{4\gamma}e^{-x}+\log 2+x
        +O\left(e^{-2x}\right),
\end{aligned}
\end{equation}
where the constant terms cancel. Imposing $Y=a_0(x)+b_0(x)Z$, with $a_0(x)=x$ and $b_0(x)=1$ in \eqref{eq:log_conditional_cont} removes the linear terms in $x$.

We now try and remove the next-order terms by using $a_1(x)=(1+\varepsilon)x$, $\varepsilon=\varepsilon(x)=o(1)$, and $b_1(x)=1+\delta$, $\delta=\delta(x)=o(1)$. Starting from \eqref{eq:log_conditional_cont}, we get
\begin{equation}\label{eq:log_conditional_final}
\begin{aligned}
  &-\dfrac{1}{2}e^{-x}\left\{
          1+\dfrac{1}{4\gamma}e^{-x}
          +e^{-\{\varepsilon x+(1+\delta)Z\}/\gamma}
            \left(
              1+\dfrac{1}{4\gamma}e^{-\{(1+\varepsilon)x+(1+\delta)Z\}}
            \right)\right\}^{\gamma}\\
  &-(1-\gamma)\log\left\{
          1+\dfrac{1}{4\gamma}e^{-x}
          +e^{-\{\varepsilon x+(1+\delta)Z\}/\gamma}
          \left(
              1+\dfrac{1}{4\gamma}e^{-\{(1+\varepsilon)x+(1+\delta)Z\}}
          \right)\right\}+\dfrac{1+\gamma}{4\gamma}e^{-x} + O\left(e^{-2x}\right)\\
=& -\dfrac{1}{2}e^{-x}\left(
              1+e^{-\{\varepsilon x+(1+\delta)Z\}/\gamma}
                          \right)^\gamma
                          \left\{
            1+\dfrac{e^{-x+\{\varepsilon x+(1+\delta)Z\}/\gamma}}{%
                  4\left(e^{\{\varepsilon x+(1+\delta)Z\}/\gamma}+1\right)}
            +\dfrac{e^{-\{(1+\varepsilon)x+(1+\delta)Z\}}}{%
                  4\left(e^{\{\varepsilon x+(1+\delta)Z\}/\gamma}+1\right)}
                          \right\}\\
  &-(1-\gamma)\left\{\log\left(
            1+e^{-\{\varepsilon x+(1+\delta)Z\}/\gamma}
                          \right)
            +\dfrac{e^{-x+\{\varepsilon x+(1+\delta)Z\}/\gamma}}{%
                  4\gamma\left(e^{\{\varepsilon x+(1+\delta)Z\}/\gamma}+1\right)}
            +\dfrac{e^{-\{(1+\varepsilon)x+(1+\delta)Z\}}}{%
                  4\gamma\left(e^{\{\varepsilon x+(1+\delta)Z\}/\gamma}+1\right)}
                          \right\}\\
  &+\dfrac{1+\gamma}{4\gamma}e^{-x} + O\left(e^{-2x}\right),
\end{aligned}
\end{equation}
We deduce from the leading order term $-(1-\gamma)\log(1+\exp[-\{\varepsilon x+(1+\delta)Z\}/\gamma])$ that $\varepsilon$  can only be $0$ as it cannot help remove any other term and comply with $\varepsilon=o(1)$. For $\delta$, there is no terms involving both $Z$ and $x$ in \eqref{eq:log_conditional_final}, that would need to be cancelled, thus $\delta=0$. Hence it is impossible to find a penultimate norming of the distribution of $\{Y-a_0(x)\}/b_0(x)\mid X=x$ for the logistic dependence structure.

\bibliographystyle{CUP}
\bibliography{bibliography}
\end{document}

%% file: macros.tex
\def\frac#1#2{{\textstyle{#1\over#2}}}

\def\sign{{\rm sign}}

\DeclareSymbolFont{AMSb}{U}{msb}{m}{n}
\DeclareMathSymbol{\Natural}{\mathbin}{AMSb}{"4E}
\DeclareMathSymbol{\Integer}{\mathbin}{AMSb}{"5A}
\DeclareMathSymbol{\Real}{\mathbin}{AMSb}{"52}
\DeclareMathSymbol{\Rational}{\mathbin}{AMSb}{"51}
\DeclareMathSymbol{\Imaginary}{\mathbin}{AMSb}{"49}
\DeclareMathSymbol{\Complex}{\mathbin}{AMSb}{"43} 
\def\bi{\begin{itemize}}
\def\ei{\end{itemize}}
\def\bd{\begin{description}}
\def\ed{\end{description}}
\def\ben{\begin{enumerate}}
\def\een{\end{enumerate}}

\def\bar#1{{\overline{#1}}}

\def\pr{{\rm Pr}}

\def\2to{{\ {\buildrel 2\over \longrightarrow}\ }}

\def\I1ton{{$I_1,\ldots,I_n$}}
\def\X1ton{{$X_1,\ldots,X_n$}}
\def\Y1ton{{$Y_1,\ldots,Y_n$}}
\def\Z1ton{{$Z_1,\ldots,Z_n$}}
\def\R1ton{{$R_1,\ldots,R_n$}}
\def\e1ton{{$e_1,\ldots,e_n$}}
\def\t1ton{{$t_1,\ldots,t_n$}}
\def\x1ton{{$x_1,\ldots,x_n$}}
\def\y1ton{{$y_1,\ldots,y_n$}}
\def\z1ton{{$z_1,\ldots,z_n$}}
\def\r1ton{{$r_1,\ldots,r_n$}}

\newcommand{\mmid}{\mathrel{}\middle|\mathrel{}}